\documentclass[10pt]{siamltex}
\usepackage{amsmath, amssymb, latexsym}
\usepackage{pstricks,pst-node}
\usepackage{srcltx}

\title {An Analysis of  broken $P_1$-Nonconforming  Finite Element Method For Interface Problems}

\author{
 Do Y. Kwak\thanks{Korea Advanced Institute of Science and Technology,
Daejeon, Korea 305-701. {{\tt email:kdy@kaist.ac.kr}}. This work was supported by a grant from Korea Science and Engineering Foundation.(R01-2007-000-10062-0) } \and K. T. Wee\thanks{Korea Advanced
Institute of Science and Technology, Daejeon, Korea 305-701 {{\tt
email:ktwee@kaist.ac.kr }}. }} \vspace{1.06in}

\date{}

\newtheorem{remark}{Remark}[section]

\makeatletter\@addtoreset{equation}{section}\makeatother

\def\pd#1#2{\frac{\partial #1}{\partial #2}}
\newcommand{\bd}{\partial}
\newcommand{\Grad}{\nabla}
\newcommand{\Div}{\mathrm{div}}

\newcommand{\wtilde}{\widetilde}
\newcommand{\what}{\widehat}
\newcommand{\dstyle}{\displaystyle}
\newcommand{\tnorm[1]}{|\!|\!| #1 |\!|\!|}

\begin{document}
\maketitle \maketitle \pagestyle{myheadings} \thispagestyle{plain} \markboth{DO
Y. KWAK AND KYE T. WEE} {$P_1$-Nonconforming  Finite Element Method For
Interface Problems}

\begin{abstract}
We study some numerical methods for solving second order elliptic problem with
interface. We introduce an immersed interface finite element method based on
the `broken' $P_1$-nonconforming piecewise linear polynomials on interface
triangular elements having edge averages as degrees of freedom. This linear polynomials are broken to match the homogeneous jump condition along the interface which is allowed to cut through the element. We prove
optimal orders of convergence in $H^1$ and $L^2$-norm. Next we propose a mixed
finite volume method in the context introduced in \cite{Kwak2003}  using the
Raviart-Thomas mixed finite element and this `broken' $P_1$-nonconforming
element. The advantage of this mixed finite volume method is that once we
solve the symmetric positive definite  pressure equation(without Lagrangian multiplier), the velocity can be computed locally by a simple
formula. This procedure avoids  solving the saddle point problem. Furthermore, we show
optimal error estimates  of velocity and pressure in our mixed finite volume method. Numerical results
show optimal orders of error in $L^2$-norm and broken $H^1$-norm for the
pressure, and in $H(\Div)$-norm for the velocity.
\end{abstract}


\noindent {\bf AMS(MOS) subject classifications.} 65N15, 65N30,
35J60.

\section{Introduction}

There are many physical problems where the underlying partial
differential equations have an interface. For example, second order
elliptic equations with discontinuous coefficients are often used to
model problems in material sciences and porous media when two or
more distinct materials or media with different conductivities,
densities or permeability  are involved. The solution of these
interface problems must satisfy interface jump conditions due to
conservation laws. If the interface is smooth enough, then the
solution of the interface problem is also  smooth in individual
regions where the coefficient is smooth, but due to the jump of the
coefficient across the interface, the global regularity is usually
low and the solution usually belongs to $H^{1+\alpha}(\Omega)$ for
some $ 0\leq\alpha<1$. Because of the low global regularity,
achieving accuracy is difficult with standard finite element
methods, unless the elements fit with the interface of general
shape.

Immersed interface method using uniform grids has many advantages
over usual fitted grid method. Using uniform grid, one does not need
to generate a grid. This is quite convenient in several aspects.
First of all,  the structure of stiffness matrix is the same as that
of the standard finite element method, where many known efficient
solvers can be exploited. Second, when a moving interface problem is
involved one does not need to generate a new grid as time evolves.
This saves considerable amount of time and storage.

The first attempt to avoid fitted grid for interface problem was
made by LeVeque and Z. Li \cite{Li1994}, where they proposed an
immersed interface method for finite difference method where the
jump condition was properly incorporated in the scheme. Cartesian
grids is most natural in this case. They subsequently applied the
same idea to other interface problems such as the Stokes flow
problem, one-dimensional moving interface problem and Hele-Shaw
flow, etc. \cite{Hele, LiStokes, Lewave, Limoving} The resulting
linear systems from these methods are non-symmetric and indefinite
even when the original problem is self-adjoint and uniformly
elliptic. Although these methods were demonstrated to be very
effective, convergence analysis of related finite difference methods
are extremely difficult and are still open.

For finite element methods, T. Lin {\em et. al.} \cite{LI:2006,
Li2004, Li2003, TLin2001} recently studied an immersed interface
finite element method using uniform grids and they proved the
approximation property of the finite element space of their scheme.
Their numerical examples demonstrated optimal orders of the error.
Other related works in this direction can be found in
\cite{Camp_LLS, Hou_L, Lai_Li_L,Lifast,Liu_SiT,Qiao_LT} and
references therein.

On the other hand, $P_1$-nonconforming finite element method introduced in
\cite{Crouzeix} for solving Stokes equation is being widely used in solving
elliptic equations  and shown to be quite effective \cite{Kwak2003, Chou2000,
Crouzeix, Kang}. Especially, it is extremely useful in solving mixed finite
element method by hybridization \cite{Arbogast, Arnold} or finite volume
formulation \cite{Kwak2003, Chou2000, ChouReport, Courbet}.

The mixed finite element method based on the dual formulation is
well-known \cite{BDFM, BDM, BrezziFortin, Chen1993, ChenDouglas,
DouglasRoberts, Falk, Girault, RT}. The motivation of mixed method
is to obtain an accurate approximation of the flow variable and
  has been widely used in the study of flow in porous
media such as petroleum engineering, underground water flow, and
electrodynamics, etc. (e.g., \cite{ChenHuan, Ewing, Pietra}) But
this scheme leads to a saddle point problem for which many well
known fast iterative methods fails.   To overcome this difficulty,
mixed hybrid methods have been introduced \cite{Arnold, Chen1993,
Marini}, where the problem reduces to a symmetric positive definite
system in Lagrange multiplier only. The flow and pressure variables are obtained via  some post -processing.

Recently, there has been some development of the mixed finite element method in
another direction: A mixed finite volume method was proposed in \cite{Courbet}
and extended in \cite{Kwak2003, Chou2000}. In this method, one use   Raviart-Thomas space and $P_1$-nonconforming space as trial spaces for velocity and pressure, and integrates the  mixed system of equation on each volume.  Then one can eliminate
velocity variable  and obtain the equation of pressure variable only (in terms of
$P_1$-nonconforming FEM) directly from the formulation without using Lagrange
multiplier. The resulting linear system is  again symmetric positive definite, and velocity can be
recovered from pressure locally in a simple manner.

The purpose of this paper is two-folded. First, we propose a
 finite element method on a uniform
triangular grid using `broken' $P_1$-functions having degrees of freedom on edges.  This is a Galerkin type $P_1$-nonconforming finite element
method  with  the basis  functions having the average  on edges as degrees of freedom, broken along the interface to match the flux condition. Then we show optimal error estimates in $H^1, L^2$-norms.
Here, we
emphasize that the meaning of `nonconforming' is different from the context of
Li et al. \cite{Li2003,LI:2006} where the basis function  has degrees of
freedom at vertices, discontinuous along edges of interface elements. Meanwhile, the basis function  here are Crouzeix-Raviart type\cite {Crouzeix}. Hence it is discontinuous along all edges intrinsically. Furthermore,
 since we use the average of  linear function(possibly broken) along edges as degrees
of freedom, the overhead of dealing with nonconformity in the proof of error estimate is significantly reduced.(See section 3)

Next, we propose a mixed finite volume method using Raviart-Thomas
space and the immersed interface finite element introduced above.
This is similar to the scheme studied in \cite{Kwak2003}, but the
usual nonconforming basis function is replaced by a `broken' one on
the interface element. We provide an optimal error analysis of
pressure and velocity.

The rest of the paper is organized as follows. In the next section, we will
describe the model problem and some preliminaries. We construct an immersed
interface $P_1$-nonconforming space with average degrees of freedom which
preserves flux continuity weakly along the interface, and prove an
interpolation error estimate. In Sections \ref{sec:IIM} and \ref{sec:L2}, we
propose an immersed interface finite element scheme and prove $H^1$ and
$L^2$-error estimates. In Section \ref{sec:mfvm}, we propose a mixed finite
volume method using Raviart-Thomas mixed finite element and our
$P_1$-nonconforming immersed interface finite element method, where the problem
reduces to symmetric positive definite system in pressure variables. The
velocity can be computed locally after pressure computation. Finally, in
Section \ref{sec:numer}, some numerical results are presented which indicate
optimal orders convergence of our methods.

\section{Preliminaries}
\begin{figure}[t]
\begin{center}
      \psset{unit=2cm}
      \begin{pspicture}(-1,-1)(1,1)
        \pspolygon(1,1)(-1,1)(-1,-1)(1,-1)
        \pscurve(0.5,0)(0.3,0.3)(-0.6,-0.1)(-0.2,-0.5)(0.52,-0.08)(0.5,0)
        \rput(0,0){\scriptsize$\Omega^-$}
        \rput(-0.3,0.5){\scriptsize$\Omega^+$}
        \rput(0.75,0){\scriptsize$\Gamma$}
      \end{pspicture}
\caption{A sketch of the domain $\Omega$ for the interface problem}
\label{fig:domain}
\end{center}
\end{figure}
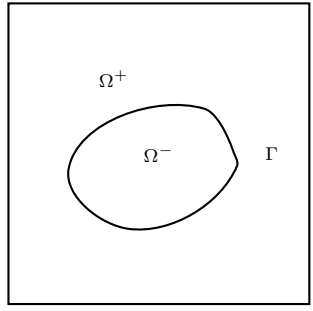

Let $\Omega$ be a convex polygonal domain in $\mathbb{R}^2$ which is separated
into two sub-domains $\Omega^+ $ and $\Omega^-$ by a $C^2$-interface $\Gamma =
\partial \Omega^- \subset\Omega,$ with $\Omega^+ = \Omega \setminus \Omega^-$
as in Fig. \ref{fig:domain}. We consider the following elliptic interface
problem
\begin{eqnarray} \label{modelpb}
\left\{
  \begin{array}{rl}
    -\Div(\beta\Grad p) = f & \hbox{ in $\Omega\setminus\Gamma$,} \\
    p = 0 & \hbox{ on $\bd\Omega$}
  \end{array}
\right.
\end{eqnarray}
with the jump conditions on the interface
\begin{eqnarray} \label{flux}
[p]=0, ~~~ [\,\beta\frac{\bd p}{\bd n}\,]=0 ~~ \textrm{ across
}\Gamma,
\end{eqnarray}
where $f\in L^2(\Omega)$ and $p\in H^1_0 (\Omega)$. We assume that
the coefficient $\beta$ is positive and piecewise constant, that is,
$ \beta(x) = \beta^- ~~ \text{for } x\in\Omega^-;~~ \beta(x) =
\beta^+ ~~ \text{for } x\in\Omega^+.$

We take as usual the weak formulation of the interface problem: Find $p\in
H^1_0(\Omega)$ such that
\begin{equation} \label{op}
\int_\Omega \beta \nabla p\cdot \nabla q dx = \int_\Omega f q dx ,
~~~ \forall q \in H^1_0(\Omega).
\end{equation}
Now we introduce the space
\begin{eqnarray*}
\wtilde{H}^2(\Omega) := \{p\in H^1(\Omega)\, :\, p\in H^2(\Omega^s),
s=+,- \}
\end{eqnarray*}
equipped with the norm
\begin{eqnarray*}
\|p\|^2_{\wtilde{H}^2(\Omega)} :=
\|p\|^2_{H^1(\Omega)}+\|p\|^2_{H^2(\Omega^+)} +
\|p\|^2_{H^2(\Omega^-)},~~ \forall\, p\in\wtilde{H}^2(\Omega),
\end{eqnarray*}
where $H^m(\Omega^s) = W^m_2(\Omega^s)$ is the usual Sobolev space of order
$m$. By Sobolev embedding theorem, for any $p\in\wtilde{H}^2(\Omega)$, we have
$p\in W^1_s(\Omega),\,\forall s>2$. Then we have the following regularity
theorem for the weak solution $p$ of the variational problem (\ref{op});
see \cite{Bramble} and \cite{Lady}.\\

\begin{theorem} \label{thm:reg}
The variational problem (\ref{op}) has a unique solution
$p\in\wtilde{H}^2(\Omega)$ which satisfies for some constant $C>0$
\begin{eqnarray}
\|p\|_{\wtilde{H}^2(\Omega)} \leq C \|f\|_{L^2(\Omega)}.
\end{eqnarray}
\end{theorem}

We now describe an immersed interface finite element method with piecewise
$P_1$-nonconforming functions.

For the simplicity of presentation, we assume that $\Omega$ is a
rectangular domain. First we consider uniform rectangular partitions
of mesh size $h$. Then we obtain triangular partitions
$\mathcal{T}_h$ by cutting the elements along diagonals. Thus we
allow the interface $\Gamma$ to cut through the elements. We assume the
following situation: The  interface
\begin{itemize}
\item   meets the edges of an interface element at no
more than two points.
\item  meets each edge  at most once except possibly it passes through two vertices.
\end{itemize}
These assumptions are reasonable if we choose $h$ sufficiently small.

  We call an element $T\in\mathcal{T}_h$ an
interface element if the interface $\Gamma$ passes through the
interior of $T$, otherwise we call $T$ a non-interface element. (If
one of the edges is part of the interface, then the element is a
non-interface element.) Let $\mathcal{E}_h$ be a collection of all
edges of $T_h$.

Let $\overline{DE}$ be the line segment connecting the intersections of the
interface and the edges of a triangle $T$. This line segment divides $T$ into
two parts $T^+$ and $T^-$ with $T=T^+\cup T^-\cup \overline{DE}$. There is a
small region in T such that $ T_r = T - (\Omega^+\cap T^+) - (\Omega^-\cap
T^-)$(see figure \ref{fig:interel}). Since $\overline{DE}$ can be considered as an approximation of the
$C^2$-curve $\Gamma\cap T$, the interface is perturbed by a $O(h^2)$ term. From
\cite{Bramble,Chen1998}, one can see for the interpolation polynomial defined
below, such a perturbation will only affect the interpolation error to the order of
$h^2$.

As usual, we want to construct local basis functions on each element $T$ of the
partition $\mathcal{T}_h$. For a non-interface element $T\in\mathcal{T}_h$, we
simply use the standard linear shape functions on $T$ having degrees of freedom
at the mid-points of the edges, and use $\overline{S}_h(T)$ to denote the
linear spaces spanned by the three nodal basis functions on $T$: Let
$m_i,\,i=1,2,3$ be the midpoints of edges of $T$. Then
\begin{eqnarray}
\overline{S}_h(T) = \text{span}\{\,\phi_i\,:\,\phi_i \text{ is linear on } T
\text{ and }\phi_i(m_j)=\delta_{ij},\,i,j=1,2,3\,\}
\end{eqnarray}
Alternatively, we can use average values along edges $e_j$ of $T$ as
degrees of freedom, i.e., $\phi_i$ can defined by
$\frac{1}{|e_j|}\int_{e_j} \phi_i\,ds=\delta_{ij},\,i,j=1,2,3$.

For this space, we have the following well-known approximation property
\cite{Ciarlet, Crouzeix}:
\begin{eqnarray} \label{noninter}
\|p - I_hp\|_{L^2(T)} + h \|p - I_hp\|_{H^1(T)} \leq C h^2 \|p\|_{H^2(T)},
\end{eqnarray}
where $I_h:H^2(T)\rightarrow \overline{S}_h(T)$ is the interpolation operator.
Finally, we use $\overline{S}_h(\Omega)$ to denote the space of the standard
piecewise $P_1$-nonconforming space with vanishing boundary nodal values.

\subsection{Local basis functions on an interface element}

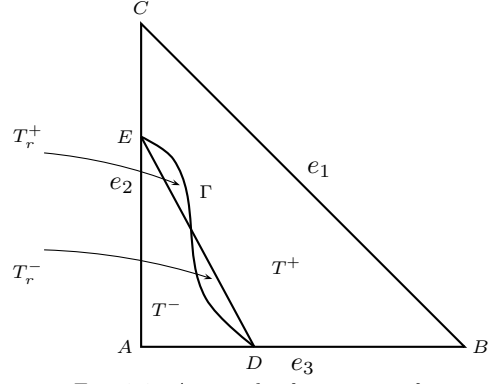
\begin{figure}[ht]
  \begin{center}
    \psset{unit=4.3cm}
    \begin{pspicture}(0,0)(1,1)
      \psset{linecolor=black} \pspolygon(0,0)(1,0)(0,1) \psline(0,0.65)(0.35,0)
      \pscurve(0,0.65)(0.1,0.58)(0.2,0.15)(0.35,0)
      \rput(0,1.05){\scriptsize$C$}
      \rput(-0.05,0){\scriptsize$A$}
      \rput(1.05,0){\scriptsize$B$}
      \rput(0.5,-0.06){$e_3$}
      \rput(0.55,0.55){$e_1$}
      \rput(-0.06,0.5){$e_2$}
      \rput(-0.05,0.65){\scriptsize$E$}
      \rput(0.08,0.12){\scriptsize$T^-$}
      \rput(0.45,0.25){\scriptsize$T^+$}
      \rput(0.35,-0.05){\scriptsize$D$}
      \rput(0.20,0.48){\scriptsize$\Gamma$}

 \rput(-.35,0.65){\scriptsize$T_r^+$}
\pnode(-.3,0.6){a}
\pnode(0.12,0.5){b}
\ncarc[linewidth=0.5\pslinewidth]{->}{a}{b}  

\pnode(-.3,.3){a}
\pnode(0.22,0.21){b}
\ncarc[linewidth=0.5\pslinewidth]{->}{a}{b}
\rput(-.35,0.23){\scriptsize$T_r^-$}
    \end{pspicture}
    \caption{A typical reference interface triangle} \label{fig:interel}
\end{center}
\end{figure}

We now consider a typical reference interface element $T$ whose geometric
configuration is given in Fig. \ref{fig:interel} in which the curve between
points $D$ and $E$ is part of the interface. Let $e_i,\,i=1,2,3$ be the edges
of $T$. For $\phi\in H^1(T)$, let $\bar{\phi}_{e_i}$ denote the average of
$\phi$ along $e_i$, i.e., $\bar{\phi}_{e_i} := \frac{1}{|e_i|}\int_{e_i}
\phi\,ds$.

We construct a piecewise linear function of the form
\begin{eqnarray}
 &&\phi(X) = \left\{%
\begin{array}{ll}
    \phi^+(X)=a_0 + b_0 x + c_0 y, & \textrm{ $X = (x,y)\in T^+$,}\\
    \phi^-(X)=a_1 + b_1 x + c_1 y, & \textrm{ $X = (x,y)\in T^-$,} \\
\end{array}%
\right. \label{def:basis-1}
\end{eqnarray}  satisfying
\begin{eqnarray}
\dstyle &&\bar{\phi}_{e_i} = V_i,~~i=1,2,3, \label{def:basis-2}\\
\dstyle &&\phi^+(D) = \phi^-(D),~\phi^+(E)=\phi^-(E), ~\beta^+ \dfrac{\bd
\phi^+}{\bd \mathbf{n}_{\overline{\textrm{\tiny{$DE$}}}}} = \beta^- \dfrac{\bd
\phi^-}{\bd
  \mathbf{n}_{\overline{\textrm{\tiny{$DE$}}}}}\, ,\label{def:basis-3}
\end{eqnarray}
where $V_i,\,i=1,2,3$ are given values and
$\mathbf{n}_{\overline{\textrm{\tiny{$DE$}}}}$ is the unit normal vector on the
line segment $\overline{DE}$. This is a piecewise linear function on $T$ that
satisfies the homogeneous jump conditions along $\overline{DE}$.

Suppose that a typical reference interface element $T$ has vertices at
$A(0,0),\,B(1,0),\,C(0,1)$. We assume that the interface meets with the edges
at $D(x_0,0)$ and $E(0,y_0)$ where $0< x_0,y_0 \leq 1$. Then the unit normal vector
to the interface is $\mathbf{n}_{\overline{\textrm{\tiny{$DE$}}}} =
(y_0,x_0)/\sqrt{x^2_0+y^2_0}$.\\

\begin{theorem} \label{thm:basis}
Given an reference interface triangle, the piecewise linear function
$\phi(x,y)$ defined by (\ref{def:basis-1})-(\ref{def:basis-3}) is
uniquely determined by three conditions $$ \bar{\phi}_{e_i}=V_i,
 \quad i=1,2,3.$$
\end{theorem}
\begin{proof}
Let $X = (x,y)^T\in T$. Since $\phi^+$ and $\phi^-$ are linear functions, we
have
\begin{equation}
\phi(X)=
\left\{%
\begin{array}{ll}
   \dstyle \phi^+(X)=a_0 + b_0 x + c_0 y, & \textrm{ $X \in T^+$,} \\
   \dstyle \phi^-(X)=a_1 + b_1 x + c_1 y, & \textrm{ $X \in T^-$.} \\
\end{array} \label{eq:basis}
\right.
\end{equation}
The condition (\ref{def:basis-2}) gives the following three equations:
\begin{eqnarray} \label{eq:basis1}
\bar{\phi}_{e_1} &=& a_0 + \frac 12 b_0 + \frac 12 c_0 = V_1
\end{eqnarray}
and
\begin{eqnarray} \nonumber
\bar{\phi}_{e_2}&=&\int_{e_2}\phi\,ds = \int_{\overline{AE}}\phi^-\,ds + \int_{\overline{EC}}\phi^+\,ds\\
&=&(a_1 + \frac{y_0}{2}c_1)y_0 + (a_0 + \frac{y_0+1}{2}c_0)(1-y_0) =
V_2, \label{eq:basis2}
\end{eqnarray}
where we used mid-point quadrature on $\overline{AE}$ and
$\overline{EC}$. Similarly, we have
\begin{eqnarray} \label{eq:basis3}
\bar{\phi}_{e_3}&=&(a_1 + \frac{x_0}{2}b_1)x_0 + (a_0 +
\frac{x_0+1}{2}b_0)(1-x_0) = V_3.
\end{eqnarray}
From the continuity condition at $D$ and $E$, we have
\begin{eqnarray}
a_0 + b_0 x_0 &=& a_1 + b_1 x_0,  \label{eq:basis4}\\
a_0 + c_0 y_0 &=& a_1 + c_1 y_0 \label{eq:basis5}
\end{eqnarray}
and the flux continuity condition along $\overline{DE}$ gives
\begin{eqnarray}
(b_0,c_0)\cdot(y_0,x_0) &=& \rho(b_1,c_1)\cdot(y_0,x_0),
\label{eq:basis6}
\end{eqnarray}
where $\rho = \beta^- / \beta^+$ and we have used that the normal
direction of the line segment $\overline{DE}$ is $(y_0,x_0)$.

Then the coefficient matrix of the above linear system for the
unknowns $a_0,b_0,c_0$ and $a_1,b_1,c_1 $ in this order is
\begin{equation}
\mathcal{A} = \left(%
\begin{array}{cccccc}
  1 & \frac 12 & \frac 12 & 0 & 0 & 0 \\
  1-y_0 & 0 & \frac12 (1-y^2_0) & y_0 & 0 & \frac12 y^2_0 \\
  1-x_0 & \frac12 (1-x^2_0) & 0 & x_0 & \frac12 x^2_0 & 0 \\
  -1 & -x_0 & 0 & 1 & x_0 & 0 \\
  -1 & 0 & -y_0 & 1 & 0 & y_0 \\
  0 & -y_0 & -x_0 & 0 & \rho y_0 & \rho x_0 \\
\end{array}%
\right)
\end{equation}
 Tedious calculation(see appendix) shows that the determinant of the matrix is
\begin{equation}\label{eq:determ}
 \textit{det}(\mathcal{A}) =\frac14(x_0^2+y_0^2)\{\rho(x_0y_0-1)-x_0y_0\}<0.
\end{equation}
  Thus the coefficients of (\ref{eq:basis})
are uniquely determined.
\end{proof}\\\

\begin{remark}
If $\bar{\phi}_{e_1},\bar{\phi}_{e_2}$ and $\bar{\phi}_{e_3}$ have
the same value, then the piecewise linear function $\phi$
satisfying (\ref{def:basis-2})-(\ref{def:basis-3}) reduces to a
constant by uniqueness.
\end{remark}\\

Now we can construct nodal basis functions on an interface element
$T$ in general position through affine mapping. We let
$\what{S}_h(T)$ to denote the three-dimensional linear space spanned
by these shape functions. We note that $\what{S}_h(T)$ is a subspace
of $H^1(T)$. Finally, we define the {\em immersed interface finite
element space} $\what{S}_h(\Omega)$ as the collection of functions
such that
\begin{displaymath}
\left\{\begin{array}{l}
\dstyle ~~\phi|_{T}\in \overline{S}_h(T),~\text{if $T$ is a noninterface element,}\\
\dstyle ~~\phi|_{T}\in \what{S}_h(T),~\text{if $T$ is an interface element,}\\
\dstyle ~~\int_e \phi|_{T_1}\,ds = \int_e \phi|_{T_2}\,ds,~\text{if $T_1,\,T_2$
are adjacent elements and $e$ is a
common edge of $T_1$ and $T_2$,}\\
\dstyle ~~\int_e \phi\,ds = 0,~\text{if $e$ is part of the boundary
$\bd\Omega$.}
\end{array}\right.
\end{displaymath}

Although for functions in $\what{S}_h(T)$ the flux jump condition is
enforced on line segments, they actually satisfy a weak flux jump
condition along the interface. This is stated in the following lemma
\cite{Li2004}, whose proof is a simple application of the divergence
theorem.\\

\begin{lemma} \label{thm:wflux}
For an interface triangle T, every function $\phi\in \what{S}_h(T)$ satisfies
the flux jump condition on $\Gamma\cap T$ in the following weak sense:
\begin{eqnarray}
\int_{\Gamma\cap T} (\beta^- \Grad \phi^- - \beta^+ \Grad \phi^+)\cdot
\mathbf{n}_{\Gamma} ds = 0.
\end{eqnarray}
\end{lemma}
\begin{proof}
Let $\phi$ be any function in $\what{S}_h(T)$. By the divergence theorem, we
have
\begin{eqnarray*}
&&\int_{\Gamma\cap T} {(\beta^- \Grad \phi^- - \beta^+ \Grad \phi^+)\cdot
\mathbf{n}_{\Gamma}} ds + \int_{\overline{DE}} {(\beta^- \Grad \phi^- - \beta^+
\Grad \phi^+)\cdot
\mathbf{n}_{\overline{\textrm{\tiny{$DE$}}}}} ds \\
&=& \int_{T_r} \Div(\beta^- \Grad \phi^- - \beta^+ \Grad \phi^+) \,dx = 0.
\end{eqnarray*}
By the flux continuity of $\phi$ on $\overline{DE}$,
\begin{eqnarray*}
\int_{\overline{DE}} {(\beta^- \Grad \phi^- - \beta^+ \Grad \phi^+)\cdot
\mathbf{n}_{\overline{\textrm{\tiny{$DE$}}}}} ds = 0,
\end{eqnarray*}
which completes the proof.
\end{proof}\\
%
\subsection{Approximation property of nonconforming immersed
interface space $\what{S}_h(T)$}

\begin{figure}[ht]
  \begin{center}
    \psset{unit=4cm}
    \begin{pspicture}(-1.5,-0.17)(1.5,1)
      \psset{linecolor=black} \pspolygon(-1.25,0)(-0.25,0)(-1.25,1)
      \psline(-1.25,0.65)(-0.9,0)
      \rput(-1.3,0){\scriptsize$A$}
      \rput(-0.2,0){\scriptsize$B$}
      \rput(-1.25,1.05){\scriptsize$C$}
      \rput(-1.12,0.12){\scriptsize$T^-$}
      \rput(-0.8,0.25){\scriptsize$T^+$}
      \rput(-0.9,-0.05){\scriptsize$D$}
      \rput(-1.3,0.65){\scriptsize$E$}
      \psset{linecolor=black} \pspolygon(0.25,0)(1.25,0)(0.25,1)
      \pscurve(0.25,0.65)(0.35,0.58)(0.45,0.15)(0.6,0)
      \rput(0.2,0){\scriptsize$A$}
      \rput(1.3,0){\scriptsize$B$}
      \rput(0.25,1.05){\scriptsize$C$}
      \rput(0.365,0.1){\scriptsize$T\cap\Omega^-$}
      \rput(0.7,0.25){\scriptsize$T\cap\Omega^+$}
      \rput(0.6,-0.05){\scriptsize$D$}
      \rput(0.2,0.65){\scriptsize$E$}
      \rput(0.45,0.48){\scriptsize$\Gamma$}
      \rput(-0.8,-0.16){\scriptsize{(a) $\what{S}_h(T)\subset H^2(T^+)\cap H^2(T^-)$}}
      \rput(0.8,-0.16){\scriptsize{(b) $\wtilde{H}^2(T)\subset H^2(T\cap\Omega^+)\cap H^2(T\cap\Omega^-)$}}
    \end{pspicture}
    \caption{Regions of $H^2$-regularity} \label{fig:regul}
\end{center}
\end{figure}
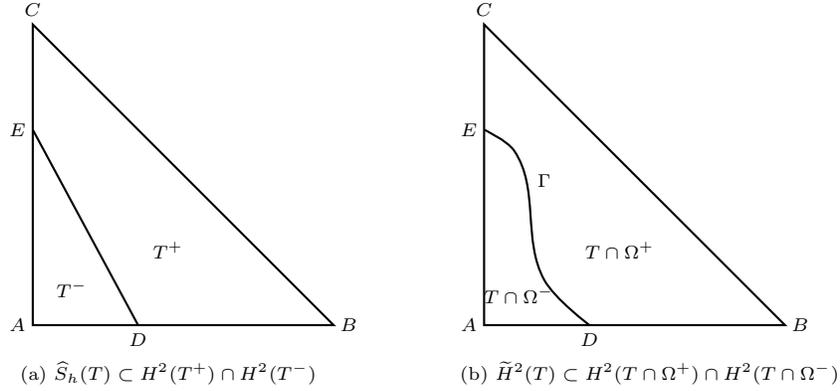

In this subsection, we would like to study the approximation property of
$\what{S}_h(T)$ by defining an interpolation operator. The difficulty lies in
the fact that $\what{S}_h(T)$ does not belong to $\wtilde{H}^2(T)$, the
restriction of $\wtilde{H}^2(\Omega)$ on $T$, where $\wtilde{H}^2(T)=H^1(T)\cap
H^2(T\cap\Omega^+)\cap H^2(T\cap\Omega^-)$(see Fig. \ref{fig:regul}). To
overcome the difficulty, we introduce a bigger space which contains both of
these spaces.

For a given interface element T, we consider a function space $X(T)$ such that
every $p\in X(T)$ satisfies
\begin{displaymath}
\left\{
\begin{array}{ll}
    \dstyle p\in H^1(T) \cap H^2(T^+\cap\Omega^+) \cap H^2(T^-\cap\Omega^-) \cap H^2(T^+_r) \cap H^2(T^-_r), \\
    \dstyle \int_{\Gamma\cap T}(\beta^- \Grad p^- -\beta^+ \Grad p^+)\cdot \mathbf{n}_\Gamma\, ds = 0, \\
\end{array}
\right.
\end{displaymath}
where $T^s_r = T_r\cap\Omega^s,~s=+,-$. For any $p\in X(T)$, we define the
following norms.
\begin{eqnarray*}
|p|^2_{X(T)} &=& |p|^2_{2,T^-\cap\Omega^-} + |p|^2_{2,T^+\cap\Omega^+} +
|p|^2_{2,T^-_r} + |p|^2_{2,T^+_r}, \\
\|p\|^2_{X(T)} &=& \|p\|^2_{1,T} + |p|^2_{X(T)}, \\
\tnorm[p]_{2,T} &=& |p|_{X(T)} + \sum^3_{i=1}{|\bar{p}_{e_i}|},
\end{eqnarray*}
where $\bar{p}_{e_i},\, i=1,2,3$ are the average on each edge $e_i$.

\begin{remark}\label{remark:X(T)}
If $p\in\wtilde{H}^2(\Omega)$, then $p|_T \in X(T)$ and
$|p|_{X(T)}=|p|_{\wtilde{H}^2(T)}$, where
$|p|^2_{\wtilde{H}^2(T)}=|p|^2_{2,T\cap\Omega^-} + |p|^2_{2,T\cap\Omega^+}$.
\end{remark}\\\

\begin{lemma} \label{lem:norm}
$\tnorm[\cdot]_{2,T}$ is a norm in the space $X(T)$ which is equivalent to
$\|\cdot\|_{X(T)}$.
\end{lemma}
\begin{proof}
Let $p\in X(T)$. If $\tnorm[p]_{2,T} = 0$, then $|p|_{X(T)} = 0$ and $|\bar{p}_{e_i}| = 0, ~ i=1,2,3$. Hence  $p$ is linear on each of the
four regions $T^+\cap\Omega^+$, $T^-\cap\Omega^-$, $T^+_r$ and $T^0_r$.  
 Since $p\in H^1(T)$, $p$ is linear on each $T^s,\,s=+,-$. Since
$p$ satisfies flux continuity condition, $\what{S}_h(T)$. Now we apply Theorem \ref{thm:basis} to
conclude $p = 0$.

We now show the equivalence of $\tnorm[\cdot]_{2,T}$ and $\|\cdot\|_{X(T)}$
(cf. \cite[p.77]{Braess}, \cite{TLin2001}). First, note that by Sobolev
embedding theorem, $H^2(T_i)$ is compactly embedded in $W^1_s(T_i)$
for any
$s>2$, where $T_i\subset T$. So we see that $X(T)\subset W^1_s(T)\subset
C^0(T)$. If $p\in X(T)$, then $p$ is a continuous function on $T$ and
$|\bar{p}_{e_i}|\leq C\|p\|_{X(T)},~i=1,2,3, $ thus
\begin{eqnarray}
\tnorm[p]_{2,T}&\leq& C\|p\|_{X(T)}, \label{eq:norms}
\end{eqnarray}
where $C$ is independent of $p$.

Now suppose that the converse
\begin{eqnarray*}
\|p\|_{X(T)}\leq C\tnorm[p]_{2,T},~~\forall p\in X(T)
\end{eqnarray*}
fails for any $C>0$. Then there exists a sequence $\{p_k\}$ in $X(T)$ with
\begin{eqnarray}\label{eq:Kondra1}
\|p_k\|_{X(T)}=1,~~\tnorm[p_k]_{2,T}\leq \frac{1}{k},~~k=1,2,\cdots .
\end{eqnarray}
Since $W^1_s(T),\,(s>2)$ is compactly imbedded in $H^1(T)$  by Kondrasov theorem
\cite[p. 114]{Ciarlet}, there exists a subsequence of $\{p_k\}$
which converges in $H^1(T)$. Without loss of generality, we can
assume that the sequence itself converges. Then $\{p_k\}$ is a
Cauchy sequence in $H^1(T)$. Noting that $|p_k|_{X(T)}\rightarrow 0$
and $\|p_k -p_l\|^2_{X(T)}\leq \|p_k -p_l\|^2_{1,T} +
(|p_k|_{X(T)}+|p_l|_{X(T)})^2$, we see that $\{p_k\}$ is a Cauchy
sequence in $X(T)$. By completeness, it converges to an element
$p^*\in X(T)$, and (\ref{eq:norms}),(\ref{eq:Kondra1}) gives
 $$ \tnorm [p^*]_{2,T}\le \tnorm [p^*-p_k]_{2,T}+\tnorm [p_k]_{2,T}\le  C\|{p^*-p_k}\|_{X(T)}+\frac1k\to 0.$$ But
\begin{eqnarray*}
\|p^*\|_{X(T)}=1~~\text{and}~~\tnorm[p^*]_{2,T}=0.
\end{eqnarray*}
This is a contradiction, since $\tnorm[p^*]_{2,T}=0$ implies $p^*=0$.
\end{proof}\\\

For any $p\in X(T)$, we define $I_h p\in \what{S}_h(T)$ using the average of
$p$ on each edge by
$$ (\overline{I_h p})_{e_i} = \bar{p}_{e_i},~~ i=1,2,3$$
and call $I_h p$ the \emph{interpolant} of $p$ in $\what{S}_h(T)$. We then
define $I_h p$ for $p\in\wtilde{H}^2(\Omega)$ by $(I_h p)|_{T} = I_h (p|_T)$.\\

\begin{lemma}
Let $T$ be an interface element. Then for any $p\in X(T)$, we have
\begin{eqnarray}
\|p - I_h p\|_{m,T} \leq C h^{2-m} \|p\|_{X(T)}, ~~ m = 0,1,
\end{eqnarray}
where $h$ is the mesh size of $T$.
\end{lemma}
\begin{proof}
Let $\what{T}$ be a reference interface element. Then for any $\hat{p}\in
X(\what{T})$
\begin{eqnarray*}
\tnorm[ \hat{p} - I_h \hat{p} ]_{2,\widehat{T}} &=&  | \hat{p} - I_h \hat{p}
|_{X(\widehat{T})} + \sum^3_{i=1} |(\overline{\hat{p} - I_h
\hat{p}})_{e_i}| \\
&=& | \hat{p} - I_h \hat{p} |_{X(\widehat{T})} =  | \hat{p} |_{X(\widehat{T})},
\end{eqnarray*}
where we used the fact that $\bar{\hat{p}}_{e_i} = (\overline{I_h
\hat{p}})_{e_i}$ on each edge and $H^2$-seminorm of the piecewise linear
function $I_h \hat{p}$ vanishes. Applying the scaling argument for $m=0,1$, we
have
\begin{eqnarray*}
\|p - I_h p\|_{m,T} &\leq& C h^{1-m} \|\hat{p} -
I_h\hat{p}\|_{m,\widehat{T}} \leq C h^{1-m} \tnorm[
\hat{p} - I_h \hat{p} ]_{2,\widehat{T}} \\
&\leq& C h^{1-m} | \hat{p} |_{X(\widehat{T})} \leq C h^{2-m} | p |_{X(T)}.
\end{eqnarray*}
\end{proof}\\\

By  above lemma, Remark \ref{remark:X(T)} and (\ref{noninter}), we obtain
the following interpolation estimate.\\

\begin{theorem}\label{thm:apperror}
For any $p\in \wtilde{H}^2(\Omega)$, there exists a constant $C>0$ such that
\begin{eqnarray}
\|p - I_h p\|_{L^2(\Omega)} + h\|p - I_h p\|_{1,h} \leq C h^2
\|p\|_{\wtilde{H}^2(\Omega)},
\end{eqnarray}
where $\|\cdot\|^2_{1,h}:=\sum_{T\in\mathcal{T}_h} \|\cdot\|^2_{1,T}$.
\end{theorem}

\section{Immersed interface FEM with `broken' $P_1$-nonconforming
elements} \label{sec:IIM}

We are now ready to define our immersed interface finite element method based
on `broken' $P_1$-nonconforming element: Find $p_h\in\what{S}_h(\Omega)$ such
that
\begin{eqnarray} \label{dp}
a_h(p_h,\phi_h) &=& (f,\phi_h),~~\forall \phi_h\in\what{S}_h(\Omega),
\end{eqnarray}
where
\begin{eqnarray}
a_h(p,\phi) &=& \sum_{T\in \mathcal{T}_h} \int_T \beta \Grad p \cdot
\Grad \phi dx,~~~ \forall\,p,\phi\in H_h(\Omega),\\
&&H_h(\Omega):= H^1_0(\Omega)+ \what{S}_h(\Omega). \nonumber
\end{eqnarray}
Here, $H_h(\Omega)$ is endowed with the piecewise $H^1$-norm $\|\cdot\|_{1,h}$.
Note that if discrete Poincar\'{e} inequality holds, then noting that the
bilinear operator $a_h(\cdot,\cdot)$ is bounded and coercive on
$\what{S}_h(\Omega)$,
the discrete problem (\ref{dp}) has a unique solution $p_h\in\what{S}_h(\Omega)$.\\

\begin{lemma}[Discrete Poincar\'{e} inequality]
There exists a constant $C>0$ independent of $h$ such that for any
$\phi\in\what{S}_h(\Omega)$
\begin{eqnarray}
C\|\phi\|^2_{L^2(\Omega)}\leq a_h(\phi,\phi).
\end{eqnarray}
\end{lemma}
\begin{proof}
Let $e$ be the common edge of two adjacent elements $T_1$ and $T_2$. Note that
since $\int_e \phi_1\,ds = \int_e \phi_2\,ds$, where
$\phi_i=\phi|_{T_i},\,i=1,2$, there exists a point $x_0\in e$ such that $\phi_1
(x_0) = \phi_2 (x_0)$. Then a slight modification of Lemma 2.1 in
\cite{Kwak1998} proves the inequality.
\end{proof}\\\

For the energy-norm error estimate of the immersed interface finite element
method, we need the well-known second Strang Lemma which is valid since
$a_h(\cdot,\cdot)$ is coercive.\\

\begin{lemma}[Second Strang Lemma] \label{lem:Strang}
If $p \in \wtilde{H}^2(\Omega),\,p_h \in \what{S}_h(\Omega) $ are the solutions
of (\ref{op}) and (\ref{dp}) respectively, then there exists a constant $C>0$
such that
\begin{eqnarray}
\| p - p_h \|_{1,h} \leq C \left\{ \inf_{q_h\in \what{S}_h(\Omega)} \| p-q_h
\|_{1,h} + \sup_{\phi_h\in \what{S}_h(\Omega)}
\frac{|\,a_h(p,\phi_h)-(f,\phi_h)\,|}{\|\phi_h\|_{1,h} }\, \right\}.
\end{eqnarray}
\end{lemma}\\\

We shall need the following estimate; see Lemma 3 in \cite{Crouzeix}.\\

\begin{lemma}\label{lem:edge-h1}
Let $e$ be an edge of $T$. Then there exists a constant $C>0$ such that for all
$\phi,\,v\in H^1(T)$
\begin{eqnarray*}
\left| \int_e \phi(v-\overline{v}_e)\,ds \right| &\leq& C h |\phi|_{1,T}
|v|_{1,T},
\end{eqnarray*}
where $\overline{v}_e := \frac{1}{|e|}\int_e v \,ds$.
\end{lemma}\\\
\begin{remark}\label{rema31}
This lemma also  holds when $\phi$ belongs to $H^1(T^s), s=\pm$ with $|\phi|_{1,T}$ understood as
sum of piecewise norm $|\phi|_{1,T^\pm}$.
\end{remark}
\begin{theorem}
Let $p\in \wtilde{H}^2(\Omega),~p_h\in \what{S}_h(\Omega)$ be the solutions of
(\ref{op}) and (\ref{dp}) respectively. Then there exists a constant $C>0$ such
that
\begin{eqnarray}
\|p-p_h\|_{1,h}\leq C h \|p\|_{\wtilde{H}^2(\Omega)}.
\label{result1}
\end{eqnarray}
\end{theorem}
\begin{proof}
We use the second Strang Lemma. The first term is nothing but an approximation
error. By Theorem \ref{thm:apperror}, we have
\begin{equation}
\inf_{q_h\in \what{S}_h(\Omega)} \| p- q_h \|_{1,h} \leq C h
\|p\|_{\wtilde{H}^2(\Omega)}.
\end{equation}
For the consistency error, we have from the definition of $a_h(\cdot,\cdot)$
and Green's formula
\begin{eqnarray}
a_h(p,\phi_h) - (f,\phi_h) &=& \sum_{T\in \mathcal{T}_h}\int_{T} {\beta \Grad
p\cdot\Grad \phi_h}\,dx - \int_{\Omega} f \phi_h \,dx \nonumber
\\ &=& \sum_{T\in \mathcal{T}_h}\int_{T} {\beta \Grad p\cdot \Grad \phi_h}\,dx - (\sum_{T\in
  \mathcal{T}_h}\int_{T} {\beta \Grad p \cdot \Grad \phi_h}\,dx - \sum_{T\in
  \mathcal{T}_h}<\beta \pd pn, \phi_h>_{\bd T})\nonumber \\ &=& \sum_{T\in
  \mathcal{T}_h}<\beta \pd pn, \phi_h>_{\bd T},
\end{eqnarray}
where $\phi_h\in \what{S}_h(\Omega)$ and $n$ is a unit outward normal
vector on each $\bd T$.
Since $\beta \pd pn$ belongs to $H^1(T\cap\Omega^+)\cap H^1(T\cap\Omega^-)$  and
 $\phi_h\in\what{S}_h(\Omega)$ has well-defined average value on
the interior edges, and vanishing average on the boundary,we have  by Lemma 
\ref{lem:edge-h1} and remark \ref{rema31}
\begin{eqnarray}
\sum_{T\in \mathcal{T}_h} <\beta \pd pn, \phi_h>_{\bd T} &=&
\sum_{T\in\mathcal{T}_h}\sum_{e\subset\bd T} <\beta \pd pn - (\overline{\beta \pd pn})_e, \phi_h>_e \nonumber\\
&\leq& \sum_{T\in\mathcal{T}_h} C h |\beta \pd pn|_{1,T} |\phi_h|_{1,T} \nonumber\\
&\leq& C h \|p\|_{\wtilde{H}^2(\Omega)} \|\phi_h\|_{1,h}. \label{4.7}
\end{eqnarray}
This completes the proof.
\end{proof}\\\


\section{$L^2$-error estimate}\label{sec:L2}

We now apply the duality argument to obtain $L^2$-norm estimate of the error.
Let us consider an auxiliary problem: Given $g\in L^2(\Omega)$, find
$\varphi\in\wtilde{H}^2(\Omega)$ such that
\begin{eqnarray}
-\Div(\beta\Grad \varphi) &=& g ~~ \textrm{ in } \Omega\setminus\Gamma, \\
\varphi &=& 0 ~~ \textrm{ on } \bd \Omega \nonumber
\end{eqnarray}
with jump conditions $[u]=0, ~~ [\,\beta\pd un\,]=0$ across $\Gamma$. Then we
have
\begin{eqnarray}
\|\varphi\|_{\wtilde{H}^2(\Omega)}\le C \|g\|_{L^2(\Omega)}.
\end{eqnarray}

Let $\varphi_h\in\what{S}_h(\Omega)$ be the solution of the corresponding
variational problem
\begin{eqnarray}
a_h(v_h,\varphi_h) = (v_h,g),~~~ \forall v_h\in\what{S}_h(\Omega).
\end{eqnarray}
Then
\begin{eqnarray*}
(p-p_h,g) &=& \sum_{T\in\mathcal{T}_h}\int_T
  \beta\Grad(p-p_h)\cdot\Grad\varphi \,dx -
  \sum_{T\in\mathcal{T}_h}\int_{\bd T} (p-p_h)\beta\pd {\varphi}{n}\,ds
  \\ &=& a_h(p-p_h,\varphi-\varphi_h) +
  a_h(p-p_h,\varphi_h) - \sum_{T\in\mathcal{T}_h}\int_{\bd T}
  (p-p_h)\beta\pd{\varphi}{n}\,ds \\ &=&
  a_h(p-p_h,\varphi-\varphi_h) + \sum_{T\in\mathcal{T}_h} \int_{\bd
    T} \beta\pd pn \varphi_h\,ds - \sum_{T\in\mathcal{T}_h}\int_{\bd T}
  (p-p_h)\beta\pd{\varphi}{n}\,ds \\ &=:& I + II - III.
\end{eqnarray*}
By continuity of $a_h(\cdot,\cdot)$ and $H^1$-error estimate of
$\varphi-\hat\varphi_h$,
\begin{eqnarray*}
|I| &\leq& C\|p-p_h\|_{1,h}\|\varphi-\varphi_h\|_{1,h}
\leq C h \|p-p_h\|_{1,h}\|\varphi\|_{\wtilde{H}^2(\Omega)} \\
&\leq& C h^2 \|p\|_{\wtilde{H}^2(\Omega)} \|\varphi\|_{\wtilde{H}^2(\Omega)}.
\end{eqnarray*}
Applying the analysis for the consistency error of $H^1$-error
estimate (\ref{4.7}), we get
\begin{eqnarray*}
|II| &=& \left|\sum_{T\in\mathcal{T}_h}\int_{\bd T} \beta\pd pn \,\varphi_h
\,ds\right| =
\left|\sum_{T\in\mathcal{T}_h}\int_{\bd T} \beta\pd pn (\varphi_h - \varphi) \,ds\right| \\
&\leq& C h \|p\|_{\wtilde{H}^2(\Omega)} \| \varphi_h-\varphi \|_{1,h} \leq C
h^2 \|p\|_{\wtilde{H}^2(\Omega)} \|\varphi\|_{\wtilde{H}^2(\Omega)}
\end{eqnarray*}
and
\begin{eqnarray*}
|III| &=& \left|\sum_{T\in\mathcal{T}_h}\int_{\bd T}
(p-p_h)\beta\pd{\varphi}{n}\,ds\right| \leq C h \|p - p_h \|_{1,h}
\|\varphi\|_{\wtilde{H}^2(\Omega)}\\ &\leq& C h^2 \|p\|_{\wtilde{H}^2(\Omega)}
\|\varphi\|_{\wtilde{H}^2(\Omega)}.
\end{eqnarray*}
Since $\|\varphi\|_{\wtilde{H}^2(\Omega)}\leq C\|g\|_{L^2(\Omega)}$, we see
that
\begin{eqnarray}
\|p-p_h\|_{L^2(\Omega)} = \sup_{g\in L^2(\Omega)}
\frac{(p-p_h,g)}{\|g\|_{L^2(\Omega)}} \leq C h^2 \|p\|_{\wtilde{H}^2(\Omega)}.
\end{eqnarray}
Thus we obtain the following $L^2$-error estimate.\\

\begin{theorem} Let $p\in \wtilde{H}^2(\Omega), ~ p_h \in
  \what{S}_h(\Omega)$ be the solutions of (\ref{op}) and (\ref{dp}) respectively.
  Then there exists a constant $C>0$ such that
\begin{equation}
\|p-p_h\|_{L^2(\Omega)} \leq C h^2 \|p\|_{\wtilde{H}^2(\Omega)}.
\end{equation}
\end{theorem}

\section{Mixed finite volume method based on IIFEM} \label{sec:mfvm}

In this section, we propose a new mixed finite volume method based on the `broken'
$P_1$-nonconforming interface finite element method introduced in the previous section.
Our method is similar to the mixed finite volume method studied in
\cite{Kwak2003, Chou2000, Courbet}, but the usual nonconforming finite element
space is replaced by our `broken' $P_1$-nonconforming space.

Let us write the problem (\ref{op}) in a mixed form by introducing
the vector variable $\mathbf{u}=-\beta\Grad p$ as
\begin{eqnarray}\label{mixedpb}
 \left\{
  \begin{array}{rl}
\dstyle    \mathbf{u}+\beta\Grad p = 0 & \hbox{ in $\Omega$,} \\
\dstyle    \Div \mathbf{u} = f & \hbox{ in $\Omega$,} \\
\dstyle    p = 0 & \hbox{ on $\bd\Omega$.}
  \end{array}
\right.
\end{eqnarray}
The mixed finite element method based on this dual formulation is
well-known \cite{BDM, BrezziFortin, RT}. The idea of the mixed
method is to find a direct approximation of the flow variable
$\mathbf{u}$. For that purpose, we introduce $\mathbf{V} =
\mathbf{H}(\Div,\Omega) =
\{\mathbf{v}\in\mathbf{L}^2(\Omega):~\Div\mathbf{v}\in
L^2(\Omega)\}$, and use the local $RT_0$ space to approximate the
flow variable which is given by $\mathbf{V}_h(T) =
\{v\,:\,v=(a+cx,b+cy),\,a,b,c\in\mathbb{R}\}$ for any triangle
element $T$. The global space $\mathbf{V}_h$ is defined as
\begin{equation}
\mathbf{V}_h=\{\mathbf{v}\,:\,\mathbf{v}|_T\in \mathbf{V}_h(T);\,
\mathbf{v}\cdot\mathbf{n} \mbox{ is continuous along interior
edges}\}.
\end{equation}

This method gives a good approximation of the flow variable.
However, it leads to a saddle point problem, that is, one obtains an
indefinite matrix system when (\ref{mixedpb}) is discretized. As
mentioned earlier, a popular way to avoid this indefinite system is
to use Lagrange multipliers\cite{Arnold}. Another possibility is to form a mixed
finite volume method as in \cite{Kwak2003, Chou2000, Courbet}.

To define a mixed finite volume method for an interface problem, we use the
well-known $RT_0$ space $\mathbf{V}_h$ for velocity and `broken'
$P_1$-nonconforming immersed interface space $\what{S}_h$ for pressure
variable. Note that every $\mathbf{v}\in\mathbf{V}_h$ has continuous normal
components across the edges of $\mathcal{T}_h$, which are constant.

We consider the following scheme: Find
$(\mathbf{u}_h,p_h)\in\mathbf{V}_h\times \what{S}_h$ which satisfies
on each element $T\in\mathcal{T}_h$
\begin{eqnarray}\label{mfvm}
\left\{
\begin{array}{rl}
\dstyle\int_T (\mathbf{u}_h+\beta\Grad p_h)\cdot\Grad\phi = 0, &
\forall\phi\in\what{S}_h(T), \\
\dstyle\int_T \Div \mathbf{u}_h = \int_T f.
\end{array}
\right.
\end{eqnarray}
Note that since $\Div\mathbf{u}_h$ is constant,
$\Div\mathbf{u}_h|_T=\overline{f}|_T:= \frac{1}{|T|}\int_T f$, where $|T|$
denotes the area of $T$. When the interface is not present, $\what{S}_h(T) =
\overline{S}_h(T)$ and this scheme coincides with the one in \cite{Kwak2003,
Courbet}. Since the numbers of unknowns and equations do not change, our scheme
is a square linear system and has a unique solution. We refer to \cite{Kwak2003} for details.

Now since $\mathbf{u}_h\cdot\mathbf{n}$ is constant on the edge   and
$\phi\in\what{S}_h$ has common average values on interior edges and vanishing
boundary nodal values, we obtain
\begin{eqnarray}
\sum_{T\in\mathcal{T}_h}\int_T \mathbf{u}_h\cdot\Grad\phi &=&
\sum_{T\in\mathcal{T}_h}\left[\int_{\bd T}(\mathbf{u}_h\cdot\mathbf{n})\phi -
\int_T \Div \mathbf{u}_h\phi\right] = - \int_\Omega \overline{f}\phi,
\end{eqnarray}
where $\overline{f}\in L^2(\Omega)$ is a simple function having value
$\overline{f}|_T$ for each $T$. From (\ref{mfvm}), it immediately follows that
\begin{eqnarray}\label{simpleform}
\sum_{T\in\mathcal{T}_h}\int_T \beta\Grad p_h\cdot\Grad\phi &=& \int_\Omega
\overline{f}\phi,~~\forall\phi\in\what{S}_h(\Omega).
\end{eqnarray}
This is the interface finite element method introduced in the previous
section, except that on the right-hand side $f$ is replaced by $\overline{f}$.

The velocity $\mathbf{u}_h$ can be computed directly from the solution $p_h$ of
(\ref{simpleform}) as follows. Let $T$ be an any element of $\mathcal{T}_h$
with the edges $e_i,~i=1,2,3$ and let $\phi_i\in\what{S}_h(T)$ be the `broken'
$P_1$-nonconforming basis function associated with the edge $e_i$. Then the
flux through the edge $e_i$ is given by
\begin{eqnarray*}
|e_i|(\mathbf{u}_h\cdot\mathbf{n})|_{e_i} &=& \int_{\bd T}
(\mathbf{u}_h\cdot\mathbf{n})\phi_i = \int_T
\Div(\mathbf{u}_h\phi_i) \\
&=& \int_T (\Div\mathbf{u}_h\phi_i + \mathbf{u}_h\cdot\Grad\phi_i),
\end{eqnarray*}
where $\phi_i\in\what{S}_h(T)$ is a basis function on $T$. Then it follows by
(\ref{mfvm}) that
\begin{eqnarray}
|e_i|(\mathbf{u}_h\cdot\mathbf{n})|_{e_i} &=& \int_T \overline{f}
\phi_i - \int_T \beta\Grad p_h\cdot\Grad\phi_i .
\end{eqnarray}
Thus in order to compute the fluxes through the edges of an element
$T$, we only need to compute the local residual of the solution
$p_h$ on each $T$.

The error estimate of $\mathbf{u}_h$ would follow that of $p_h$. In
fact, we can relate the estimate $\|\mathbf{u}-\mathbf{u}_h\|_0$
with $\|p-p_h\|_{1,h}$. First, we show the following local formula.\\

\begin{lemma}
Let $\mathbf{u}_h,\, p_h$ be the solutions of (\ref{mfvm}), then
\begin{eqnarray}\label{formula}
\mathbf{u}_h(\mathbf{x})&=&-\overline{\beta\Grad
p_h}+\frac{\overline{f}}{2} (\mathbf{x}-\mathbf{x}_B),
~~\forall\mathbf{x}\in T,
\end{eqnarray}
where $\overline{\beta\Grad p_h}$ denotes the average of $\beta\Grad
p_h$ on $T$ and $\mathbf{x}_B$ is the center of $T$.
\end{lemma}
\begin{proof}
Expanding $\mathbf{u}_h$ about $\mathbf{x}_B$, the barycenter of
$T$, we have
\begin{eqnarray} \label{taylor}
\mathbf{u}_h(\mathbf{x}) &=& \mathbf{u}_h(\mathbf{x}_B) +
\mathcal{D}\mathbf{u}_h(\mathbf{x}_B)(\mathbf{x}-\mathbf{x}_B),~~\mathbf{x}\in
T,
\end{eqnarray}
where $\mathcal{D}\mathbf{u}_h$ is the Jacobian matrix of $\mathbf{u}_h$. Let
$\mathbf{u}_h=(a+cx,b+cy)\in \mathbf{V}_h(T)$, then we have
\begin{eqnarray}
\mathcal{D}\mathbf{u}_h(\mathbf{x}_B)(\mathbf{x}-\mathbf{x}_B) &=&
c(\mathbf{x}-\mathbf{x}_B) \nonumber =
\frac{\Div\mathbf{u}_h}{2}(\mathbf{x}-\mathbf{x}_B) =
\frac{\overline{f}}{2}(\mathbf{x}-\mathbf{x}_B),
\end{eqnarray}
where we used the relation $\Div\mathbf{u}_h=\overline{f}$. On the other hand,
applying $\phi=(x,0)^T$ and $(0,y)^T$ in (\ref{mfvm}), we see
\begin{eqnarray}
-\overline{\beta\Grad p_h} &=& \frac{1}{|T|}\int_T \mathbf{u}_h =
\mathbf{u}_h(\mathbf{x}_B).
\end{eqnarray}
Substituting these into (\ref{taylor}), we obtain formula
(\ref{formula}).
\end{proof}\\\

\begin{remark}
Our formula is different from the one in \cite{Kwak2003, Chou2000} where they have
\begin{eqnarray}\label{formula0}
\mathbf{u}_h(\mathbf{x})&=&- {\beta\Grad
p_h}+\frac{\overline{f}}{2} (\mathbf{x}-\mathbf{x}_B),
~~\forall\mathbf{x}\in T.
\end{eqnarray}
The $p_h$ in
our scheme is broken along a line segment contained in an interface element, Hence we have taken the average of $\beta\Grad p_h$.
\end{remark}\\\

By the above lemma, we have
\begin{eqnarray*}
\mathbf{u}(\mathbf{x})|_T-\mathbf{u}_h(\mathbf{x}) &=& -\beta\Grad p
+ \overline{\beta\Grad p_h} -\frac{\overline{f}}{2}
(\mathbf{x}-\mathbf{x}_B).
\end{eqnarray*}
So
\begin{eqnarray*}
\|\mathbf{u}-\mathbf{u}_h\|_{0,T}&\leq&\|\beta\Grad p -
\overline{\beta\Grad
p_h}\|_{0,T}+C|\overline{f}|\|\mathbf{x}-\mathbf{x}_B\|_{0,T}.
\end{eqnarray*}
Since $\beta$ is piecewise constant and $\mathbf{u}=-\beta\Grad p$,
we have
\begin{eqnarray*}
\|\beta\Grad p - \overline{\beta\Grad p_h}\|_{0,T} &\leq&
\|\beta\Grad p - \overline{\beta\Grad p}\|_{0,T} +
\|\overline{\beta\Grad p - \beta\Grad p_h}\|_{0,T} \\
&\leq& C h \|\mathbf{u}\|_{1,T}+C h |\overline{\beta\Grad p - \beta\Grad p_h}|\\
&\leq& C h \|\mathbf{u}\|_{1,T}+C|p-p_h|_{1,T},
\end{eqnarray*}
provided $\mathbf{u}\in \mathbf{H}^1(T)$. Hence
\begin{eqnarray*}
\|\mathbf{u}-\mathbf{u}_h\|_{0,T}&\leq& C
\{h \|\mathbf{u}\|_{1,T}+|p-p_h|_{1,T}+ h^2 |\overline{f}|\} \\
&\leq& C\{h \|\mathbf{u}\|_{1,T}+|p-p_h|_{1,T}+h\|f\|_{0,T}\},
\end{eqnarray*}
where we used $|\overline{f}|\leq h^{-1}\|f\|_{0,T}$ and
$\|\mathbf{x}-\mathbf{x}_B\|_{0,T}\leq C h^2$.

Summing over every $T\in\mathcal{T}_h$, we have
\begin{eqnarray}
\|\mathbf{u}-\mathbf{u}_h\|_{L^2(\Omega)}&\leq& C\{h
\|\mathbf{u}\|_{\mathbf{H}^1(\Omega)}+\|p-p_h\|_{1,h}+h\|f\|_{L^2(\Omega)}\}.
\end{eqnarray}
Since $\Div\mathbf{u}=f$ and $\Div\mathbf{u}_h=\overline{f}$, we can
easily obtain the estimate for
$\|\Div\mathbf{u}-\Div\mathbf{u}_h\|_{L^2(\Omega)}$. This is summarized in the following.\\

\begin{theorem} \label{thm:errflux}
Let $\mathbf{u}_h,\, p_h$ be the solutions of (\ref{mfvm}), then there exists a
constant $C>0$ such that
\begin{eqnarray}
\|\mathbf{u}-\mathbf{u}_h\|_{L^2(\Omega)}+\|\Div \mathbf{u}-\Div
\mathbf{u}_h\|_{L^2(\Omega)} &\leq& C h \{
\|\mathbf{u}\|_{\mathbf{H}^1(\Omega)}+\|p\|_{\wtilde{H}^2(\Omega)}+
\|f\|_{1,h}\},
\end{eqnarray}
provided $\mathbf{u}\in \mathbf{H}^1(\Omega)$.
\end{theorem}\\\

\section{Numerical examples} \label{sec:numer}
In this section, we report numerical results for the schemes
introduced previously. For a numerical test, we solve problem
(\ref{modelpb}) with the rectangular domain $\Omega =
[-1,1]\times[-1,1]$ partitioned into unform right triangles having
step size $h$. We take a circle with radius $r_0 = 0.5$ as an
interface, and the exact solution is chosen as
\begin{eqnarray*}
 p = \left\{%
\begin{array}{ll}
    \dfrac{r^3}{\beta^-} & \textrm{in $\Omega^-$,} \\
    \dfrac{r^3}{\beta^+}+(\dfrac{1}{\beta^-}-\dfrac{1}{\beta^+})r^3_0 & \textrm{in $\Omega^+$.} \\
\end{array}%
\right.
\end{eqnarray*}

We note that this example is taken from Z. Li \cite{Li2004}. We present errors
in $L^2, H^1$-norm for the pressure $p$, while in $H(\Div)$-norm for the
velocity $\mathbf{u}$. Here the order of convergence is determined by the least
squares fit to the data. In Table \ref{table:1-1000} and \ref{table:1000-1}(first two columns), we report
the results of the `broken' $P_1$-nonconforming immersed interface scheme
introduced in Section \ref{sec:IIM}, where we used conjugate gradient
method(CG) to solve the resulting discrete system. It shows   optimal
order  of convergence for $L^2$-norm and $H^1$-norm:
\begin{eqnarray}
\|p-p_h\|_0 &\approx& O(h^2),~~~\|p-p_h\|_{1,h} \approx O(h).
\end{eqnarray}
We also present some result for mixed finite volume method introduced in
Section \ref{sec:mfvm}. Again, this shows optimal order of convergence for the
flow variable which is consistent with Theorem \ref{thm:errflux} (cf. last columns of Table
\ref{table:1-1000} and \ref{table:1000-1}):
\begin{eqnarray}
\|\mathbf{u}-\mathbf{u}_h\|_{L^2(\Omega)} + \|\Div (\mathbf{u}-
\mathbf{u}_h\|_{L^2(\Omega)} \approx O(h).
\end{eqnarray}
This is in good agreement with some fitted grid computation; when the jump of
the coefficient is large, one usually have $O(h)$ order accuracy, see
\cite{Cai} problem 1, p.310, for example.

\begin{table}[hbt]
\begin{center}
\begin{tabular}{|c|c|c|c|c|c|c|c|c|}
\hline $1/h$ & $\|p-p_h\|_{0}$ & order & $\|p-p_h\|_{1,h}$ & order &
 $\|\mathbf{u}-\mathbf{u}_h\|_0$ & order & $\|\Div(\mathbf{u}-\mathbf{u}_h)\|_0$ & order \\
\hline 8  & 9.576e-3 &       & 1.208e-1 &       &  2.945e-1 &        & 1.053e+0 &  \\
\hline 16 & 2.666e-3 & 1.845 & 6.744e-2 & 0.841 &  1.702e-1 & 0.791  & 5.292e-1 & 0.993 \\
\hline 32 & 6.488e-4 & 2.039 & 3.341e-2 & 1.013 &  8.906e-2 & 0.934  & 2.650e-1 & 0.998 \\
\hline 64 & 1.400e-4 & 2.212 & 1.657e-2 & 1.012 &  4.290e-2 & 1.054  & 1.326e-1 & 0.999 \\
\hline 128 & 3.716e-5 & 1.914 & 8.242e-3 & 1.008 &  2.015e-2 & 1.090 & 6.629e-2 & 1.000 \\
\hline 256 & 8.973e-6 & 2.050 & 4.117e-3 & 1.001 &  9.865e-3 & 1.030 & 3.315e-2 & 1.000 \\
\hline Order & & 2.029 & &  0.985 & &   0.994 & & 0.998 \\
\hline
\end{tabular}
\caption{Nonconforming immersed interface FEM: $\beta^-=1,\beta^+
=1000$} \label{table:1-1000}
\end{center}
\end{table}

\begin{table}[hbt]
\begin{center}
\begin{tabular}{|c|c|c|c|c|c|c|c|c|}
\hline $1/h$ & $\|p-p_h\|_{0}$ & order & $\|p-p_h\|_{1,h}$ & order &
 $\|\mathbf{u}-\mathbf{u}_h\|_0$ & order & $\|\Div(\mathbf{u}-\mathbf{u}_h)\|_0$ & order\\
\hline 8  & 1.447e-2 &       & 6.575e-1 &       &  3.361e-1 &       & 1.053e+0 &\\
\hline 16 & 3.497e-3 & 2.049 & 3.312e-1 & 0.989 &  1.657e-1 & 1.020 & 5.292e-1 & 0.993\\
\hline 32 & 8.826e-4 & 1.986 & 1.661e-1 & 0.996 &  8.165e-2 & 1.021 & 2.650e-1 & 0.998\\
\hline 64 & 2.210e-4 & 1.998 & 8.311e-2 & 0.999 &  4.075e-2 & 1.003 & 1.326e-1 & 0.999\\
\hline 128 & 5.507e-5 & 2.005 & 4.157e-2 & 0.999 &  1.959e-2 & 1.057 & 6.629e-2 & 1.000\\
\hline 256 & 1.370e-5 & 2.007 & 2.079e-2 & 1.000 &  9.658e-3 & 1.020 & 3.315e-2 & 1.000\\
\hline Order & &  2.005 & &  0.997 & &  1.024 & & 0.998\\
\hline
\end{tabular}
\caption{Nonconforming immersed interface FEM: $\beta^-=1000,\beta^+
=1$} \label{table:1000-1}
\end{center}
\end{table}
\newpage
\appendix{\bf Appendix: Computation of determinant of ${\cal A}$. }

By adding the last three columns to first the  three, we obtain
\begin{eqnarray*}
& \left|%
\begin{array}{cccccc}
  1 & \frac 12 & \frac 12 & 0 & 0 & 0 \\
  1-y & 0 & \frac12 (1-y^2) & y & 0 & \frac12 y^2 \\
  1-x  & \frac12 (1-x^2) & 0 & x & \frac12 x^2 & 0 \\
  -1 & -x & 0 & 1 & x & 0 \\
  -1 & 0 & -y & 1 & 0 & y \\
  0 & -y & -x & 0 & \rho y & \rho x \\
\end{array}\right|
 = \left|%
\begin{array}{cccccc}
  1 & \frac 12 & \frac 12 & 0 & 0 & 0 \\
  1  & 0 & \frac12  & y & 0 & \frac12 y^2 \\
  1  & \frac12   & 0 & x & \frac12 x^2 & 0 \\
  0 &  0 & 0 & 1 & x & 0 \\
  0 & 0 &  0 & 1 & 0 & y \\
  0 & (\rho-1)y & (\rho-1)x & 0 & \rho y & \rho x \\
\end{array}\right|\end{eqnarray*}
Now gaussian elimination gives
\begin{eqnarray*}
&=\left|%
\begin{array}{cccccc}
  1 & \frac 12 & \frac 12 & 0 & 0 & 0 \\
  0 & -\frac12 &       0  & y & 0 & \frac12 y^2 \\
  0 & 0        & -\frac12 & x  & \frac12 x^2 & 0 \\
  0 &  0       & 0       & 1 & x & 0 \\
  0 & 0        &  0      & 1 & 0 & y \\
  0 &   (\rho-1)y & (\rho-1)x & 0& \rho y & \rho x \\
\end{array}\right|
 =\left| %
\begin{array}{cccccc}
  1 & \frac 12 & \frac 12 & 0 & 0 & 0 \\
  0 & -\frac12 &       0  & y & 0 & \frac12 y^2 \\
  0 & 0        & -\frac12 & x  & \frac12 x^2 & 0 \\
  0 &  0       & 0       & 1 & x & 0 \\
  0 & 0        &  0      & 1 & 0 & y \\
  0 & 0        &  (\rho-1)x &2( \rho-1)y^2 & \rho y & C \\
\end{array}\right|
\end{eqnarray*}
where $C=\rho x+ (\rho-1)y^3.$ Continuing
\begin{eqnarray*}
&=\left|%
\begin{array}{cccccc}
  1 & \frac 12 & \frac 12 & 0 & 0 & 0 \\
  0 & -\frac12 &       0  & y & 0 & \frac12 y^2 \\
  0 & 0        & -\frac12 & x  & \frac12 x^2 & 0 \\
  0 &  0       & 0       & 1 & x & 0 \\
  0 & 0        &  0      & 1 & 0 & y \\
  0 & 0        &  0 &A & B & C \\
\end{array}\right|
\end{eqnarray*}
where
$$A=2(\rho-1)(x^2+y^2),\quad B= \rho y+ (\rho-1)x^3.$$ Hence the determinant is
\begin{eqnarray*}
\frac14( xy A - y B - C x) &=&\frac14 \{2 (\rho-1) xy(x^2+y^2)-y (\rho y+
(\rho-1)x^3)-x(\rho x+ (\rho-1)y^3)\} \\
&=&\frac14 (x^2+y^2)\{\rho(xy-1)-xy\}.
 \end{eqnarray*}
This verifies (\ref{eq:determ}).

\end{document}